\documentclass[11pt]{article}
\usepackage{graphicx, colordvi}
\usepackage{dsfont}
\usepackage{epsfig}
\usepackage{mathrsfs}
\usepackage{amssymb,amsfonts,amsmath,amsthm,cite,color}
\usepackage{dsfont}
\usepackage{CJKutf8}
\usepackage{epsfig}
\usepackage{mathrsfs}
\usepackage{longtable}
\usepackage{enumerate}
\usepackage{geometry}
\usepackage{hyperref}
\hypersetup{
	colorlinks=true,
	linkcolor=cyan,
	filecolor=blue,
	urlcolor=red,
	citecolor=green,
}

\newtheorem{theo}{Theorem}

\newtheorem{lem}[theo]{Lemma}

\newtheorem{coj}[theo]{Conjecture}
\makeatletter \@addtoreset{equation}{section}
\@addtoreset{theo}{section} \makeatother

\newcommand{\bN} { {\mathbb{N}}}

\newcommand{\bZ} { {\mathbb{Z}}}

\newcommand{\qbinom}[2]{\left[ \begin{matrix} #1 \\ #2 \end{matrix} \right]}
\def\qed{\hfill \rule{4pt}{7pt}}

\author{}
\title{}

\begin{document}
\begin{center}
 {\large $q$-Congruences for Z.-W. Sun's generalized polynomials $w^{(\alpha)}_k(x)$}
\end{center}

\begin{center}
 {Lin-Yue Li}$^{1}$ and {Rong-Hua Wang}$^{2}$

   $^{1,2}$School of Mathematical Sciences\\
   Tiangong University \\
   Tianjin 300387, P.R. China\\
   $^{1}$2330141409@tiangong.edu.cn\\
   $^{2}$wangronghua@tiangong.edu.cn \\[10pt]
\end{center}

\vskip 6mm \noindent {\bf Abstract.}
In 2022, Z.-W. Sun defined
\begin{equation*}
	w_k^{(\alpha)}{(x)}=\sum_{j=1}^{k}w(k,j)^{\alpha}x^{j-1},
 \end{equation*}
where $k,\alpha$ are positive integers and $w(k,j)=\frac{1}{j}\binom{k-1}{j-1}\binom{k+j}{j-1}$.
Let $(x)_{0}=1$ and $(x)_{n}=x(x+1)\cdots(x+n-1)$ for all $n\geq 1$.
In this paper, it is proved by $q$-congruences that for any positive integers ${\alpha,\beta, m,n,r}$, we have
\begin{equation*}
	\frac{(2,n)}{n(n+1)(n+2)}\sum_{k=1}^{n}k^r(k+1)^r(2k+1)w_{k}^{(\alpha)}(x)^{m}\in\bZ[x],
\end{equation*}
\begin{equation*}
	\frac{(2,n)}{n(n+1)(n+2)}\sum_{k=1}^{n}(-1)^{k}k^r(k+1)^r(2k+1)
	w_{k}^{(\alpha)}(x)^{m}\in\bZ[x],
\end{equation*}
and
\begin{equation*} \frac{2}{[n,n+1,\cdots,n+2\beta+1]}\sum_{k=1}^{n}(k)_{\beta}^r(k+\beta+1)_{\beta}^r(k+\beta) \prod_{i=0}^{2\beta-1}w_{k+i}^{(\alpha)}(x)^m\in\bZ[x],
\end{equation*}
where $[n,n+1,\cdots,n+2\beta+1]$ is the least common multiple of $n$, $n+1$, $\cdots$, $n+2\beta+1$.
Taking $r=\beta=1$ above will confirm some of Z.-W. Sun’s conjectures.


\section{Introduction}
Let $\bN$ and $\bZ^+$ denote the set of nonnegative and positive integers respectively.
For any $n\in\bZ^{+}$, Sun\cite{Sun2018} defined $w(n,k)$ as
\begin{equation*}
w(n,k)=\frac{1}{k}\binom{n-1}{k-1}\binom{n+k}{k-1}
=\binom{n-1}{k-1}\binom{n+k}{k}-\binom{n}{k}\binom{n+k}{k-1}\in\bZ.
\end{equation*}
Clearly $w(n,n)=\binom{2n}{n}/(n+1)$ is the Catalan number.
Sun\cite{Sun2018,Sun2022} also proved that when $m$ are positive integers and $m\leq n$,
\begin{equation*}
	\sum_{k=m}^{n}(-1)^{n-k}\binom{k-1}{m-1}w(n,k)=w(n,m)
\end{equation*}
and established the following connections between $w(n,k)$ and the Narayana  number $N(n,k)=\frac{1}{n}\binom{n}{k}\binom{n}{k-1}$,
\begin{equation*}
	w(n,k)=\sum_{j=1}^{k}\binom{n-j}{k-j}N(n,j)
\quad\text{and}\quad
	N(n,k)=\sum_{j=1}^{k}\binom{n-j}{k-j}(-1)^{k-j}w(n,j)
\end{equation*}
for any integers $n\geq k\geq 1$.
Recently, the authors \cite{Li and Wang2025} showed that
$
	\sum_{k=0}^{2b}w(n,k+1)
$
is odd for any $n,b\in\bN$.

In 2018, Sun\cite{Sun2018} defined the polynomial $w_n(x)$ as
\begin{equation*}
	w_n(x)=\sum_{k=1}^{n}w(n,k)x^{k-1}
\end{equation*}
and discovered that
\begin{equation*}
    w_n(-1-x)=(-1)^{n-1}w_n(x).
\end{equation*}
Later, Sun\cite{Sun2022} proved that, for any $n\in\bZ^{+}$, $w_n(x)$ equals the little Schr\"oder
polynomial $s_n(x)=\sum_{k=1}^{n}\frac{1}{n}\binom{n}{k}\binom{n}{k-1}x^{k-1}(x+1)^{n-k}$
and that
\begin{equation*}
	(2x+1)\sum_{k=1}^{n}k(k+1)(2k+1)(-1)^{n-k}w_k(x)^2=n(n+1)(n+2)w_n(x)w_{n+1}(x).
\end{equation*}

In 2022, Sun\cite{Sun2022} further defined
\begin{equation*}
	w_n^{(\alpha)}{(x)}:=\sum_{k=1}^{n}w(n,k)^{\alpha}x^{k-1}\quad{\alpha,n}\in\bZ^{+}
\end{equation*}
and provided many interesting conjectures involving $w_n^{(\alpha)}{(x)}$.
\begin{coj}[Conjecture 5.2 of Sun\cite{Sun2022}]\label{th:2022sun Conjecture 5.2}
	For any ${\alpha,m,n}\in\bZ^{+}$, we have
	\begin{equation}\label{eq:coj w(n,k)}
		\frac{(2,n)}{n(n+1)(n+2)}\sum_{k=1}^{n}k(k+1)(2k+1)w_{k}^{(\alpha)}(x)^{m}\in\bZ[x].
	\end{equation}
Also,
	\begin{equation}\label{eq:coj (2,m-1,n)}
	\frac{(2,m-1,n)}{n(n+1)(n+2)}\sum_{k=1}^{n}(-1)^{k}k(k+1)(2k+1)w_{k}(x)^{m}\in\bZ[x],
\end{equation}
and
	\begin{equation}\label{eq:coj 1/n(n+1)(n+2)}
	\frac{1}{n(n+1)(n+2)}\sum_{k=1}^{n}(-1)^{k}k(k+1)(2k+1)w_{k}^{(\alpha)}(x)^{m}\in\bZ[x] \quad for\ \alpha >1.
\end{equation}
\end{coj}
\begin{coj}[Conjecture 5.4 of Sun\cite{Sun2022}]\label{th:2022sun Conjecture 5.4}
	\textnormal{(i)} For any ${\alpha,m,n}\in\bZ^{+}$, we have
	\begin{equation}\label{coj:(i)}
		\frac{2(2,n)}{n(n+1)(n+2)}\sum_{k=1}^{n}k(k+1)(k+2)	(w_k^{(\alpha)}(x)w_{k+1}^{(\alpha)}(x))^m\in\bZ[x].
	\end{equation}
	\textnormal{(ii)} For any ${m,n}\in\bZ^{+}$, we have
	\begin{equation}\label{coj:(ii)}
		\frac{2(2,n)}{n(n+1)(n+2)(2x+1)^m}\sum_{k=1}^{n}k(k+1)(k+2)	(w_k(x)w_{k+1}(x))^m\in\bZ[x].
	\end{equation}
	If $n\in\bZ^{+}$ is even, then
\begin{equation}\label{coj:(iii)}
	\frac{4}{n(n+1)(n+2)(2x+1)^3}\sum_{k=1}^{n}k(k+1)(k+2)	w_k(x)w_{k+1}(x)\in\bZ[x].
\end{equation}
\end{coj}
Recall that $w_{2j}(x)/(2x+1)\in\bZ[x]$ for all $j\in\bZ^+$ which has proved by Sun\cite{Sun2018}, one can directly get \eqref{coj:(ii)} from \eqref{coj:(i)} by noting that either $k$ or $k+1$ is even. Recently, \eqref{coj:(iii)} was confirmed by Jia and Huang\cite{jia and huang2025}.

In this paper, we will focus on conjectures \eqref{eq:coj w(n,k)}--\eqref{coj:(i)} by the method of $q$-analogues which are widely used in the proof and discovery of arithmetic properties of many classical combinatorial sequences,  consult \cite{Guo2006,Guo2018,Guo2019a,Guo2019b,Guo2019c,Guo2020,Liu 2020,Ni Pan 2020,Pan and Cao2006,Pan2014,Pan and sun2012,hou2021,Sun2019} and references therein.

To proceed, some preliminaries are need.
For each integer $n$, the $q$-integer denoted by $[n]_q$, or $[n]$ for short, is defined by
\begin{equation*}
	[n]=\frac{1-q^n}{1-q}.
\end{equation*}
Obviously $\lim\limits_{q\to1}[n]=n$. For $k\in\bZ$, the $q$-binomial coefficient $\qbinom{n}{k}$ is given by
\begin{equation*}
	\qbinom{n}{k}=\frac{\prod_{1\leq j\leq k}[n-j+1]}{\prod_{1\leq j\leq k}[j]}.
\end{equation*}
 Also, we set $\qbinom{n}{0}=1$ and  $\qbinom{n}{k}=0$ if $k<0$ or $k>n$.

In this paper, we define
\begin{equation}\label{eq:q-w(n,k)}	w_k^{(\alpha)}(x;q)=\sum_{j=1}^{k}q^{\alpha\left(\binom{j+1}{2}-(k+1)(j-1)\right)}
\left(\qbinom{k-1}{j-1}\qbinom{k+j}{j}-\qbinom{k}{j}\qbinom{k+j}{j-1}\right)^{\alpha}x^{j-1}.
\end{equation}
Clearly $w_k^{(\alpha)}(x;q)$ is a $q$-analogue of $w_k^{(\alpha)}(x)$.
Then we will show
\begin{theo}\label{th1:k(k+1)(2k+1):q}
	For any ${\alpha,m,n,r}\in\bZ^{+}$, we have
	\begin{equation}\label{eq1:no -1 q-analogue}
		\sum_{k=0}^{n-1}[k(k+1)]^r[2k+1]q^{(n-1-k)(\alpha m+1)}	w_k^{(\alpha)}(x;q)^m \equiv0\pmod{[n]}.
	\end{equation}
\end{theo}
It is well-known that the cyclotomic polynomial
\begin{equation*}
	\Phi_n(q)=\prod_{\substack{k=1 \\ (k,n)=1}}^{n}(q-e^{\frac{2\pi ik}{n}})\in\bZ[q]
\end{equation*}
is irreducible in the ring $\bZ[q]$. For the alternating cases, we obtain the following fact.
\begin{theo}\label{th:(-1)k(k+1)(k+2):q}
	For any ${\alpha,m,n,r}\in\bZ^{+}$, we have
	\begin{equation}\label{eq:have -1 q-analogue}
		\sum_{k=0}^{n-1}(-1)^{k}[k(k+1)]^r_{q^2}[2k+1]q^{(n-1-k)(2\alpha m+1)}	w_k^{(\alpha)}(x;q^2)^m
	\end{equation}
	is divisible by
	\begin{equation*}
		\prod_{\substack{d\mid n \\ d>1\ is\ odd}}\Phi_d(q)\cdot\prod_{\substack{d\mid n \\ d>1\ is\ even}}\Phi_d(q^2).
	\end{equation*}
\end{theo}

Taking $q=1$ in Theorem \ref{th1:k(k+1)(2k+1):q} and Theorem \ref{th:(-1)k(k+1)(k+2):q} respectively.
We have
\[
n\mid \sum\limits_{k=1}^{n-1}\varepsilon^k k^{r}(k+1)^{r}(2k+1)w_{k}^{(\alpha)}(x)^{m}
\]
for $\varepsilon=\pm1$. Then it is straightforward to verify that
\[
(n+i)\mid \sum\limits_{k=1}^{n}\varepsilon^k k^{r}(k+1)^{r}(2k+1)w_{k}^{(\alpha)}(x)^{m},
\]
where $i=0,1$ or $2$. This leads to
\begin{theo}\label{th:2022sun Conjecture 5.2 generalize}
	For any ${\alpha,m,n,r}\in\bZ^{+}$, we have
	\begin{equation}\label{eq:th w(n,k)}
	\frac{(2,n)}{n(n+1)(n+2)}\sum_{k=1}^{n}k^r(k+1)^r(2k+1)w_{k}^{(\alpha)}(x)^{m}\in\bZ[x],
	\end{equation}
	and
	\begin{equation}\label{eq:th -w(n,k)}
\frac{(2,n)}{n(n+1)(n+2)}\sum_{k=1}^{n}(-1)^{k}k^r(k+1)^r(2k+1)
w_{k}^{(\alpha)}(x)^{m}\in\bZ[x].
	\end{equation}
\end{theo}
Clearly \eqref{eq:th w(n,k)} is a generalization of Sun's \eqref{eq:coj w(n,k)}
while \eqref{eq:th -w(n,k)} generalized \eqref{eq:coj (2,m-1,n)} and \eqref{eq:coj 1/n(n+1)(n+2)} only when $n$ are odd integers.
In the cases when the summand contains the multiplication $w_k^{(\alpha)}(x;q)w_{k+1}^{(\alpha)}(x;q)$, we have
\begin{theo}\label{th:w_k+i}
For any ${\alpha,m,n,r}\in\bZ^{+}$, we have
	\begin{equation*}
		\sum_{k=0}^{n-1}[k(k+2)]^r[2(k+1)]q^{(n-2-k)(2\alpha m+1)}	(w_k^{(\alpha)}(x;q)w_{k+1}^{(\alpha)}(x;q))^m \equiv0\pmod{[n]}.
	\end{equation*}
\end{theo}
The proof of \ref{th:w_k+i} motivates the discovery of the following generalization.
\begin{theo}\label{th:generalizedw_k+i}
	For any ${\alpha,\beta, m,n,r}\in\bZ^{+}$, we have
	\begin{equation*}		\sum_{k=0}^{n-1}[(k)_{\beta}(k+\beta+1)_{\beta}]^r[2(k+\beta)]q^{(n-2\beta-k)(2\beta\alpha m+1)}\prod_{i=0}^{2\beta-1}(w_{k+i}^{(\alpha)}(x;q))^m\equiv0\pmod{[n]}.
	\end{equation*}
\end{theo}

Taking $q=1$ above gives
\[
n\mid 2\sum\limits_{k=1}^{n-1}(k)_{\beta}^r(k+\beta+1)_{\beta}^r(k+\beta) \prod_{i=0}^{2\beta-1}w_{k+i}^{(\alpha)}(x)^m.
\]
Then it is easy to check fact
\[
(n+i)\mid 2\sum\limits_{k=1}^{n}(k)_{\beta}^r(k+\beta+1)_{\beta}^r(k+\beta) \prod_{i=0}^{2\beta-1}w_{k+i}^{(\alpha)}(x)^m
\]
for any $i$ with $0\leq i \leq 2\beta+1$ when $r\geq 1$.
This leads to
\begin{theo}
	For any ${\alpha, \beta,  m,n,r}\in\bZ^{+}$, we have
	\begin{equation*}\label{eq:general WkWk+1} \frac{2}{[n,n+1,\cdots,n+2\beta+1]}\sum_{k=1}^{n}(k)_{\beta}^r(k+\beta+1)_{\beta}^r(k+\beta) \prod_{i=0}^{2\beta-1}w_{k+i}^{(\alpha)}(x)^m\in\bZ[x].
	\end{equation*}
\end{theo}
This confirms and generalizes the conjecture \eqref{coj:(i)}.

The rest of the paper is organized as follows. In section 2 and section 3, we provide the proof of Theorem \ref{th1:k(k+1)(2k+1):q} and \ref{th:(-1)k(k+1)(k+2):q} respectively. Theorem \ref{th:w_k+i} will be discussed in  section 4.

\section{Proof of Theorem \ref{th1:k(k+1)(2k+1):q}}
It is well-known that
\begin{equation*}
	[n]=\prod_{\substack{d\mid n \\ d>1}}\Phi_d(q).
\end{equation*}
Therefore, it suffices to prove that
\begin{theo}\label{th:k(k+1)(2k+1):q}
	For any ${\alpha,m,n,r}\in\bZ^{+}$ and $d>1$ being a divisor of $n$, we have
	\begin{equation}\label{eq:no -1 q-analogue}
		\sum_{k=0}^{n-1}[k(k+1)]^r[2k+1]q^{(n-1-k)(\alpha m+1)}	w_k^{(\alpha)}(x;q)^m \equiv0\pmod{\Phi_d(q)}.
	\end{equation}
\end{theo}
To this aim, we first provide another expression of $w_k^{(\alpha)}(x;q)$.
\begin{lem}
	For any $k,\alpha\in\bZ^{+}$, we have
	\begin{equation}\label{eq:q-w(n,k)another form}
		w_k^{(\alpha)}(x;q)=\sum_{-\infty<j<\infty}(-1)^{\alpha j}q^{\alpha j^2}\left(q^{k+1}\qbinom{k-1}{j-1}\qbinom{-k-1}{j}-\qbinom{k}{j}\qbinom{-k-2}{j-1}\right)^{\alpha}x^{j-1}.
	\end{equation}
\end{lem}
\noindent\emph{Proof.}  Not that
\begin{align}\label{eq:k+j,j}
	\qbinom{k+j}{j}
	=&\frac{(1-q^{k+1})(1-q^{k+2})\cdots(1-q^{k+j})}{(1-q)(1-q^{2})\cdots(1-q^{j})}\notag\\
	=&(-1)^{j}q^{kj+\binom{j+1}{2}}\frac{(1-q^{-k-1})(1-q^{-k-2})\cdots(1-q^{-k-j})}{(1-q)(1-q^{2})\cdots(1-q^{j})}\notag\\
	=&(-1)^{j}q^{kj+\binom{j+1}{2}}\qbinom{-k-1}{j}.
\end{align}
and
\begin{align}\label{eq:k+j,j-1}
	\qbinom{k+j}{j-1}
	=&\frac{(1-q^{k+2})(1-q^{k+3})\cdots(1-q^{k+j})}{(1-q)(1-q^{2})\cdots(1-q^{j-1})}\notag\\
	=&(-1)^{j-1}q^{(k+1)(j-1)+\binom{j}{2}}\frac{(1-q^{-k-2})(1-q^{-k-3})\cdots(1-q^{-k-j})}{(1-q)(1-q^{2})\cdots(1-q^{j-1})}\notag\\
	=&(-1)^{j-1}q^{(k+1)(j-1)+\binom{j}{2}}\qbinom{-k-2}{j-1}.
\end{align}
Substituting \eqref{eq:k+j,j} and \eqref{eq:k+j,j-1} into \eqref{eq:q-w(n,k)} gives the equality \eqref{eq:q-w(n,k)another form}.
\qed

To prove Theorem \ref{th:k(k+1)(2k+1):q}, we recall the well-known $q$-Lucas Theorem \cite{BES1992}.
\begin{lem}\label{lem:Lucas}
	Suppose that $d>1$ is a positive integer and $a,b,s,t$ are integers with $0\leq b,t\leq d-1$. Then, we have
	\begin{equation*}
		\qbinom{ad+b}{sd+t}\equiv {\binom{a}{s}\qbinom{b}{t}}\pmod{\Phi_d(q)}.
	\end{equation*}
\end{lem}

In the follow, we always suppose $d>1$ is a divisor of $n$ and $n=hd$.
Writing $k=ad+b$ and $j=sd+t$ where $0\leq b,t\leq d-1$.
\begin{lem}\label{lem:w tongyu B}
	Let $w_{k}^{(\alpha)}(x;q)$ be defined as in \eqref{eq:q-w(n,k)another form}.
	For $d>2$ and $1\leq b\leq d-2$, we have
	 \begin{equation}\label{eq:w_ad+b tongyu B}
	 	w_{ad+b}^{(\alpha)}(x;q)\equiv {B_{a,b,d}^{(\alpha)}(x;q)}\pmod {\Phi_d(q)}
	 \end{equation}
 and
 \begin{equation}\label{eq:w_ad+d-1-b tongyu B}
 	w_{ad+d-b-1}^{(\alpha)}(x;q)\equiv q^{-\alpha(2b+1)}{B_{a,b,d}^{(\alpha)}(x;q)}\pmod {\Phi_d(q)},
 \end{equation}
 where
 \begin{align}\label{eq:B_a,b,d}
 	B_{a,b,d}^{(\alpha)}(x;q)
 	=\sum_{\substack{-\infty<s<\infty \\  1\leq t\leq d-1}}
 	&(-1)^{\alpha(sd+t)}q^{\alpha t^2}
 	\binom{a}{s}^{\alpha}\binom{-a-1}{s}^{\alpha}\notag\\
 	&\times\left(q^{b+1}\qbinom{b-1}{t-1}\qbinom{d-b-1}{t}
 	  +\qbinom{b}{t}\qbinom{d-b-2}{t-1}\right)^{\alpha}x^{sd+t-1}.
 \end{align}
\end{lem}
\noindent\emph{Proof.}
Let \begin{equation*}
A_{a,b,s,t,d}(q)= q^{ad+b+1}\qbinom{ad+b-1}{sd+t-1}\qbinom{-ad-b-1}{sd+t}
              +\qbinom{ad+b}{sd+t}\qbinom{-ad-b-2}{sd+t-1}.
\end{equation*}
By Lemma \ref{lem:Lucas},  one can see
\begin{align*}
	&A_{a,b,s,0,d}(q)\notag\\
	\equiv&\ q^{b+1}\binom{a}{s-1}\qbinom{b-1}{d-1}\binom{-a-1}{s}\qbinom{d-b-1}{0}+	\binom{a}{s}\qbinom{b}{0}\binom{-a-1}{s-1}\qbinom{d-b-2}{d-1}\notag\\
	=&\ 0\pmod{\Phi_d(q)}
\end{align*}
since $b-1<d-1$ and $d-b-2<d-1$.
Then, we have
\begin{align*}
	w_{ad+b}^{(\alpha)}(x;q)
	=&\sum_{\substack{-\infty<s<\infty \\  1\leq t\leq d-1}}
	(-1)^{\alpha(sd+t)}q^{\alpha(sd+t)^2}(A_{a,b,s,t,d}(q))^{\alpha}x^{sd+t-1}\notag\\
    &+\sum_{-\infty<s<\infty}
    (-1)^{\alpha sd}q^{\alpha(sd)^2}
    (A_{a,b,s,0,d}(q))^{\alpha}x^{sd-1}\notag\\
	\equiv&\ B_{a,b,d}^{(\alpha)}(x;q)
	\pmod{\Phi_d(q)}.
\end{align*}

For the proof of \eqref{eq:w_ad+d-1-b tongyu B},
note that
\begin{align}
	&q^{b+1}\qbinom{b-1}{t-1}\qbinom{d-b-1}{t}+\qbinom{b}{t}\qbinom{d-b-2}{t-1}\notag\\
\equiv & \ q^{2b+1}\left(q^{d-b}\qbinom{b}{t}\qbinom{d-b-2}{t-1}
	+\qbinom{b-1}{t-1}\qbinom{d-b-1}{t}\right)\pmod{\Phi_d(q)}\notag.
\end{align}
We have
\begin{align*}
	w_{ad+d-b-1}^{(\alpha)}(x;q)
	\equiv& \ B_{a,d-b-1,d}^{(\alpha)}(x;q)\notag\\
	=&\sum_{\substack{-\infty<s<\infty \\  1\leq t\leq d-1}}
	(-1)^{\alpha(sd+t)}q^{\alpha t^2}
	\binom{a}{s}^{\alpha}\binom{-a-1}{s}^{\alpha}\notag\\
	&\times\left(q^{d-b}\qbinom{b}{t}\qbinom{d-b-2}{t-1}
	+\qbinom{b-1}{t-1}\qbinom{d-b-1}{t}\right)^{\alpha}x^{sd+t-1}\notag\\
	\equiv&\ q^{-\alpha(2b+1)}B_{a,b,d}^{(\alpha)}(x;q)\pmod{\Phi_d(q)}.
\end{align*}
This ends the proof. \qed

\emph{Proof of Theorem \ref{th:k(k+1)(2k+1):q}.}	
For $d=2$, we have
\begin{align}\label{eq:d=2 sum (1.2)}
	&\sum_{k=0}^{n-1}[k(k+1)]^r[2k+1]q^{(n-1-k)(\alpha m+1)}
	w_k^{(\alpha)}(x;q)^m\notag\\
	\equiv&\ \sum_{\substack{0\leq a \leq h-1 \\ 0\leq b \leq 1}}[b(b+1)]^{r}q^{(-1-b)(\alpha m+1)}	w_{2a+b}^{(\alpha)}(x;q)^m\notag\\
	\equiv&\ 0\pmod{\Phi_2(q)}
\end{align}
by noting that either $b$ or $b+1$ is even.
For $d>2$,
\begin{align}\label{eq:d>2}
	&\sum_{k=0}^{n-1}[k(k+1)]^r[2k+1]q^{(n-1-k)(\alpha m+1)}
	w_k^{(\alpha)}(x;q)^m\notag\\
	\equiv&\sum_{\substack{0\leq a \leq h-1 \\ 1\leq b \leq d-2}}[(ad+b)(ad+b+1)]^r[2ad+2b+1]q^{(hd-1-ad-b)(\alpha m+1)}	w_{ad+b}^{(\alpha)}(x;q)^m\notag\\
	\equiv&\sum_{\substack{0\leq a \leq h-1 \\ 1\leq b \leq d-2}}[b(b+1)]^r[2b+1]q^{(-1-b)(\alpha m+1)}	B_{a,b,d}^{(\alpha)}(x;q)^m
	\pmod{\Phi_d(q)}
\end{align}
with the help of \eqref{eq:w_ad+b tongyu B}.

Replacing $k$ by $n-1-k$ above,  we have
\begin{align}\label{eq:d>2 substitute}
	&\sum_{k=0}^{n-1}[(n-1-k)(n-k)]^r[2n-2k-1]q^{k(\alpha m+1)}
	w_{n-1-k}^{(\alpha)}(x;q)^m\notag\\
	\equiv&\sum_{\substack{0\leq a \leq h-1 \\ 1\leq b \leq d-2}}[b(b+1)]^r[-2b-1]q^{b(\alpha m+1)-\alpha m(2b+1)}B_{h-a-1,b,d}^{(\alpha)}(x;q)^m\pmod{\Phi_d(q)}\notag\\
	=&\sum_{\substack{0\leq \tilde{a} \leq h-1 \\ 1\leq b \leq d-2}}[b(b+1)]^r[-2b-1]q^{b(\alpha m+1)-\alpha m(2b+1)}B_{\tilde{a},b,d}^{(\alpha)}(x;q)^m
\end{align}
with the help of \eqref{eq:w_ad+d-1-b tongyu B}, where $\tilde{a}=h-a-1$.
Combining  \eqref{eq:d>2} and \eqref{eq:d>2 substitute} together will arrive at
\begin{align*}
	&2\sum_{k=0}^{n-1}[k(k+1)]^r[2k+1]q^{(n-1-k)(\alpha m+1)}
	w_k^{(\alpha)}(x;q)^m\notag\\
	\equiv&\sum_{\substack{0\leq a \leq h-1 \\ 1\leq b \leq d-2}}[b(b+1)]^{r}q^{-\alpha m-b\alpha m}
	\left([2b+1]q^{-1-b}+[-2b-1]q^{b}\right)	B_{a,b,d}^{(\alpha)}(x;q)^m\notag\\
	\equiv& \ 0 \pmod{\Phi_d(q)}
\end{align*}
by noting that
\begin{equation*}
	[2b+1]q^{-1-b}+[-2b-1]q^{b}=0.
\end{equation*}
This concludes the proof.\qed

\section{Proof of Theorem \ref{th:(-1)k(k+1)(k+2):q}}
The proof of the alternating cases are more complicated. To proceed, we first recall the following fact, one can see \cite[Lemma 3.1]{Pan2014} for a proof.
\begin{lem}\label{lem:d is odd or even}
	If $d>1$ is odd, then $\Phi_d(q)$ divides $\Phi_d(q^2)$. If $d$ is even, then $\Phi_d(q^2)=\Phi_{2d}(q)$.
\end{lem}

\emph{Proof of Theorem \ref{th:(-1)k(k+1)(k+2):q}.}
We will discuss according to the parity of $d$. If $d$ is odd.
Substituting $q=q^2$ in \eqref{eq:w_ad+b tongyu B}, we get
\small\begin{align}\label{eq:d>2(odd)}
	&\sum_{k=0}^{n-1}(-1)^{k}[k(k+1)]^r_{q^2}[2k+1]q^{(n-1-k)(2\alpha m+1)}
	w_k^{(\alpha)}(x;q^2)^m\notag\\
	\equiv&\sum_{\substack{0\leq a \leq h-1 \\ 1\leq b \leq d-2}}(-1)^{ad+b}[b(b+1)]^r_{q^2}[2ad+2b+1]q^{(h-a)d-(1+b)(2\alpha m+1)}	B_{a,b,d}^{(\alpha)}(x;q^2)^m\pmod{\Phi_d(q^2)}\notag\\
	\equiv&\sum_{\substack{0\leq a \leq h-1 \\ 1\leq b \leq d-2}}(-1)^{ad+b}[b(b+1)]^r_{q^2}[2b+1]q^{(-1-b)(2\alpha m+1)}	B_{a,b,d}^{(\alpha)}(x;q^2)^m\pmod{\Phi_d(q)}
\end{align}
with help of Lemma \ref{lem:d is odd or even}.
Replacing $k$ by $n-1-k$, we have
\begin{align}\label{eq:d>2substitute(odd)}
	&\sum_{k=0}^{n-1}(-1)^{n-1-k}[(n-1-k)(n-k)]^r_{q^2}[2n-2k-1]q^{k(2\alpha m+1)}
	w_{n-1-k}^{(\alpha)}(x;q^2)^m\notag\\
	\equiv&\sum_{\substack{0\leq a \leq h-1 \\ 1\leq b \leq d-2}}(-1)^{(h-a)d-1-b}[b(b+1)]^r_{q^2}[2(h-a)d-2b-1]
	 q^{ad-2b\alpha m+b-2\alpha m}B_{h-a-1,b,d}^{(\alpha)}(x;q^2)^m\notag\\
	\equiv&\sum_{\substack{0\leq \tilde{a} \leq h-1 \\ 1\leq b \leq d-2}}(-1)^{\tilde{a}d+d-1-b}[b(b+1)]^r_{q^2}[-2b-1]q^{-2b\alpha m+b-2\alpha m}B_{\tilde{a},b,d}^{(\alpha)}(x;q^2)^m\pmod{\Phi_d(q)}
\end{align}
by substituting $q=q^2$ in \eqref{eq:w_ad+d-1-b tongyu B}.

By  \eqref{eq:d>2(odd)} and \eqref{eq:d>2substitute(odd)}, we have
\begin{align*}
	&2\sum_{k=0}^{n-1}(-1)^{k}[k(k+1)]^r_{q^2}[2k+1]q^{(n-1-k)(2\alpha m+1)}	w_k^{(\alpha)}(x;q^2)^m\notag\\
	\equiv&\sum_{\substack{0\leq a \leq h-1 \\ 1\leq b \leq d-2}}(-1)^{ad+b}[b(b+1)]^r_{q^2}q^{-2\alpha m(b+1)}
	\left([2b+1]q^{-1-b}+[-2b-1]q^{b}\right)	B_{a,b,d}^{(\alpha)}(x;q^2)^m\notag\\
	=&\ 0\pmod{\Phi_d(q)}.
\end{align*}
Then we have $\Phi_d(q)$ divides \eqref{eq:have -1 q-analogue}.
Next we assume $d$ is even. For $d=2$, we have
\begin{align}\label{eq:d=2 sum (even)}
	&\sum_{k=0}^{n-1}(-1)^{k}[k(k+1)]^r_{q^2}[2k+1]q^{(n-1-k)(2\alpha m+1)}
	w_k^{(\alpha)}(x;q^2)^m\notag\\
	\equiv&\sum_{\substack{0\leq a \leq h-1 \\ 0\leq b \leq 1}}(-1)^{b}[b(b+1)]^r_{q^2}[2b+1]q^{2h-2a-(1+b)(2\alpha m+1)}w_{2a+b}^{(\alpha)}(x;q^2)^m\notag\\
	\equiv&\ 0\pmod{\Phi_2(q^2)}.
\end{align}

For $d>2$, we have
\begin{align}\label{eq:d>2(even)}
	&\sum_{k=0}^{n-1}(-1)^{k}[k(k+1)]^r_{q^2}[2k+1]q^{(n-1-k)(2\alpha m+1)}
	w_k^{(\alpha)}(x;q^2)^m\notag\\
	\equiv&\sum_{\substack{0\leq a \leq h-1 \\ 1\leq b \leq d-2}}(-1)^{ad+b}[b(b+1)]^r_{q^2}[2ad+2b+1]q^{(hd-ad)-(1+b)(2\alpha m+1)}B_{a,b,d}^{(\alpha)}(x;q^2)^m\notag\\
	\equiv&\sum_{\substack{0\leq a \leq h-1 \\ 1\leq b \leq d-2}}(-1)^{b}[b(b+1)]^r_{q^2}[2b+1]q^{(hd-ad)-(1+b)(2\alpha m+1)}B_{a,b,d}^{(\alpha)}(x;q^2)^m\pmod{\Phi_d(q^2)}
\end{align}
with the help of Lemma \ref{lem:d is odd or even}.

Similarly, we have
\begin{align}\label{eq:d>2 substitute (even)}
	&\sum_{k=0}^{n-1}(-1)^{n-1-k}[(n-1-k)(n-k)]^r_{q^2}[2n-2k-1]q^{k(2\alpha m+1)}
	w_{n-1-k}^{(\alpha)}(x;q^2)^m\notag\\
	\equiv&\sum_{\substack{0\leq a \leq h-1 \\ 1\leq b \leq d-2}}(-1)^{hd-1-ad-b}[b(b+1)]^r_{q^2}[-2b-1]q^{ad-2b\alpha m+b-2\alpha m}	B_{h-a-1,b,d}^{(\alpha)}(x;q^2)^m\pmod{\Phi_d(q^2)}\notag\\
	=&\sum_{\substack{0\leq \tilde{a} \leq h-1 \\ 1\leq b \leq d-2}}(-1)^{-1-b}[b(b+1)]^r_{q^2}[-2b-1]q^{(h-\tilde{a}-1)d-2b\alpha m+b-2\alpha m}B_{\tilde{a},b,d}^{(\alpha)}(x;q^2)^m.
\end{align}
Note that $\Phi_d(q^2)=\Phi_{2d}(q)$. So, we have $q^d\equiv-1\pmod{\Phi_d(q^2)}$.
By \eqref{eq:d>2(even)} and \eqref{eq:d>2 substitute (even)}, we have
\begin{align}\label{eq:d>2 sum even}
	&2\sum_{k=0}^{n-1}(-1)^{k}[k(k+1)]^r_{q^2}[2k+1]q^{(n-1-k)(2\alpha m+1)}	w_k^{(\alpha)}(x;q^2)^m\notag\\
	\equiv&\sum_{\substack{0\leq a \leq h-1 \\ 1\leq b \leq d-2}}(-1)^{b}[b(b+1)]^r_{q^2}q^{hd-ad-2\alpha m(b+1)}
	\left([2b+1]q^{-1-b}+[-2b-1]q^{b}\right)	B_{a,b,d}^{(\alpha)}(x;q^2)^m\notag\\
	=&\ 0\pmod{\Phi_d(q^2)}.
\end{align}
Combining \eqref{eq:d=2 sum (even)} and \eqref{eq:d>2 sum even}, we have $\Phi_d(q^2)$ divides \eqref{eq:have -1 q-analogue}.
This concludes the proof.
\qed

\section{Proof of Theorem \ref{th:w_k+i}}
By \eqref{eq:w_ad+b tongyu B} and \eqref{eq:w_ad+d-1-b tongyu B}, we can check the following congruences
\begin{equation}\label{eq:w_ad+b+1 tongyu B}
	w_{ad+b+1}^{(\alpha)}(x;q)\equiv {B_{a,b+1,d}^{(\alpha)}(x;q)}\pmod {\Phi_d(q)}
\end{equation}
and
\begin{equation}\label{eq:w_ad+d-2-b tongyu B}
	w_{ad+d-b-2}^{(\alpha)}(x;q)\equiv q^{-\alpha(2b+3)}{B_{a,b+1,d}^{(\alpha)}(x;q)}\pmod {\Phi_d(q)}
\end{equation}
hold for $d>3$ and $0\leq b\leq d-3$.

\emph{Proof of Theorem \ref{th:w_k+i}.}	
We only need to prove
\begin{equation}\label{eq:n-2Wk(x)wk+1(x)}
\sum_{k=0}^{n-2}[k(k+2)]^r[k+1]q^{(n-2-k)(2\alpha m+1)}
(w_{k}^{(\alpha)}(x;q)w_{k+1}^{(\alpha)}(x;q))^m \equiv 0 \pmod{\Phi_d(q)}
\end{equation}
by noting that $[n]$ divides the term when $k=n-1$.

When $d=2$ or $3$, \eqref{eq:n-2Wk(x)wk+1(x)} holds since one of $k$, $k+1$ or $k+2$ can be divisible by $d$.
For $d>3$, we obtain
\begin{align}\label{eq:d>4}
	&\sum_{k=0}^{n-2}[k(k+2)]^r[2(k+1)]q^{(n-2-k)(2\alpha m+1)}
	(w_k^{(\alpha)}(x;q)w_{k+1}^{(\alpha)}(x;q))^m\notag\\
	\equiv&\sum_{\substack{0\leq a \leq h-1 \\ 1\leq b \leq d-3}}[(ad+b)(ad+b+2)]^r[2(ad+b+1)]q^{(hd-2-ad-b)(2\alpha m+1)}
	(w_{ad+b}^{(\alpha)}(x;q)w_{ad+b+1}^{(\alpha)}(x;q))^m\notag\\
	\equiv&\sum_{\substack{0\leq a \leq h-1 \\ 1\leq b \leq d-3}}[b(b+2)]^r[2(b+1)]q^{(-2-b)(2\alpha m+1)}	(B_{a,b,d}^{(\alpha)}(x;q)B_{a,b+1,d}^{(\alpha)}(x;q))^m
	\pmod{\Phi_d(q)}
\end{align}
with the help of \eqref{eq:w_ad+b tongyu B} and \eqref{eq:w_ad+b+1 tongyu B}.

Replacing $k$ by $n-2-k$ above,  we have
\begin{align}\label{eq:d>4 substitute}
	&\sum_{k=0}^{n-2}[(n-2-k)(n-k)]^r[2(n-k-1)]q^{k(2\alpha m+1)}
	(w_{n-2-k}^{(\alpha)}(x;q)w_{n-1-k}^{(\alpha)}(x;q))^m\notag\\
	\equiv&\sum_{\substack{0\leq a \leq h-1 \\ 1\leq b \leq d-3}}[b(b+2)]^r[-2(b+1)]q^{b(2\alpha m+1)-4\alpha m(b+1)}(B_{h-a-1,b,d}^{(\alpha)}(x;q)B_{h-a-1,b+1,d}^{(\alpha)}(x;q))^m\notag\\
	=&\sum_{\substack{0\leq \tilde{a} \leq h-1 \\ 1\leq b \leq d-3}}[b(b+2)]^r[-2(b+1)]q^{b(2\alpha m+1)-4\alpha m(b+1)}(B_{\tilde{a},b,d}^{(\alpha)}(x;q)B_{\tilde{a},b+1,d}^{(\alpha)}(x;q))^m\pmod{\Phi_d(q)}
\end{align}
with the help of \eqref{eq:w_ad+d-1-b tongyu B} and \eqref{eq:w_ad+d-2-b tongyu B}.
Combining  \eqref{eq:d>4} and \eqref{eq:d>4 substitute} together will arrive at
\begin{align}
	&2\sum_{k=0}^{n-2}[k(k+2)]^r[2(k+1)]q^{(n-2-k)(2\alpha m+1)}
	(w_k^{(\alpha)}(x;q)w_{k+1}^{(\alpha)}(x;q))^m\notag\\
	\equiv&\sum_{\substack{0\leq a \leq h-1 \\ 1\leq b \leq d-3}}[b(b+2)]^{r}q^{-2b\alpha m-4\alpha m}
	\left([2(b+1)]q^{-b-2}+[-2(b+1)]q^{b}\right)
	(B_{a,b,d}^{(\alpha)}(x;q)B_{a,b+1,d}^{(\alpha)}(x;q))^m\notag\\
	\equiv& \ 0 \pmod{\Phi_d(q)}\notag
\end{align}
by noting that
\begin{equation*}
	[2(b+1)]q^{-b-2}+[-2(b+1)]q^{b}=0.
\end{equation*}
This ends the proof.
\qed


\begin{thebibliography}{10}		
    \bibitem{Guo2018}
    V.J.W. Guo,
    \newblock{A $q$-analogue of a Ramanujan-type supercongruence involving central binomial coefficients,}
    \newblock {\em J. Math. Anal. Appl.} 458 (2018), 590--600.

    \bibitem{Guo2019a}
    V.J.W. Guo,
    \newblock{A $q$-analogue of the (I.2) supercongruence of Van Hamme,}
    \newblock {\em Int. J. Number Theory }15
    (2019), 29--36.

    \bibitem{Guo2019b}
    V.J.W. Guo,
    \newblock{A $q$-analogue of a curious supercongruence of Guillera and Zudilin,}
    \newblock {\em J. Diﬀerence Equ. Appl.} 25 (2019), 342--350.


    \bibitem{Guo2020}
    V.J.W. Guo and M.J. Schlosser,
    \newblock{A family of $q$-hypergeometric congruences modulo the
    fourth power of a cyclotomic polynomial,}
    \newblock {\em Israel J. Math.} 240 (2020), 821--835.


	\bibitem{Guo2006}
     V.J.W. Guo and J. Zeng,
    \newblock{Some arithmetic properties of the $q$-Euler numbers and $q$-Salié numbers,}
   \newblock {\em European J. Combin.} 27 (2006), 884--895.

    \bibitem{Guo2019c}
    V.J.W. Guo and W. Zudilin,
    \newblock{A $q$-microscope for supercongruences,}
    {\em Adv. Math.} 346 (2019), 329--358.

    \bibitem{hou2021}
    Q.-H. Hou and Z.-W. Sun,
    \newblock{$q$-Analogues of some series for powers of $\pi$,}
    \newblock{\em Ann. Comb.} 25(2021), 167--177.

	\bibitem{jia and huang2025}
    C.-B. Jia and J.-Q. Huang,
    \newblock{Congruences involving Delannoy numbers and Schröder numbers,}
    \newblock{arXiv:2410.17522v1}.

    \bibitem{Liu 2020}
    J.-C. Liu,
    \newblock{On a congruence involving $q$-Catalan numbers,}
    \newblock {\em C. R. Math. Acad. Sci. Paris }358(2020), 211--215 .

	\bibitem{Li and Wang2025}
    L.-Y. Li and R.-H. Wang,
    \newblock{Arithmetic properties of generalized Delannoy polynomials and   Schröder polynomials,}
    \newblock{arXiv:2503.12748v1}.


    \bibitem{Ni Pan 2020}
    H.-X. Ni and H. Pan,
    \newblock{Some symmetric $q$-congruences modulo the square of a cyclotomic polynomial,}
    \newblock {\em J. Math. Anal. Appl. }481 (2020), Art. 123372.

	
	
	
	
	\bibitem{Pan2014}
	H. Pan,
	\newblock{On divisibility of sums of Apéry polynomials,}
	\newblock{\em J. Number Theory} 143 (2014), 214--223.

	\bibitem{Pan and Cao2006}
	H. Pan and H.-Q. Cao,
	\newblock{A congruence involving products of $q$-binomial coefficients,}
	\newblock {\em J. Number Theory} 121 (2006), 224--233.

    \bibitem{Pan and sun2012}
    H. Pan and Z.-W Sun,
    \newblock{Some $q$-congruences related to $3$-adic valuations,}
    \newblock {\em Adv. Appl. Math.} 49(2012), 263--270.


	\bibitem{BES1992}
    Bruce E. Sagan,
    \newblock{Congruence properties of $q$-analogs,}
    \newblock {\em Adv. Math.} 95 (1992), 127--143.
		
    \bibitem{Sun2018}
     Z.-W. Sun,
    \newblock{Arithmetic properties of Delannoy numbers and Schr\"oder numbers,}
    \newblock{\em J. Number Theory} 183 (2018), 146--171.

    \bibitem{Sun2019}
    Z.-W. Sun,
    \newblock{Two $q$-analogues of Euler’s formula $\zeta(2)=\pi^{2}/6$,}
   \newblock{\em Colloq. Math.} 158(2019), 313--320.

	
	\bibitem{Sun2022}
	Z.-W. Sun,
	\newblock{On Motzkin numbers and central trinomial coefficients,}
	\newblock {\em Adv. Appl. Math.} 136 (2022), Article ID 102319.
	
\end{thebibliography}
\end{document}